\documentclass[12pt,a4paper]{amsart}

\usepackage{amsmath,amssymb,amsthm,amsfonts}
\renewcommand\eqref[1]{(\ref{#1})} 

\usepackage[mathscr]{eucal}
\usepackage{enumitem} 

\setitemize[1]{label=$\bullet$,leftmargin=*}

\newcommand{\abs}[1]{\lvert#1\rvert}

\newcommand{\absBig}[1]{\Bigl\lvert#1\Bigr\rvert}

\newcommand{\set}[1]{\left\{#1\right\}}

\newcommand\R{{\mathbb R}}
\newcommand\C{{\mathbb C}}
\newcommand\N{{\mathbb N}}
\newcommand\NUO{{\mathbb N}\cup\{0\}}

\newcommand\RN{{{\mathbb R}^N}}

\newcommand\SNm{{{\mathbb S}^{N-1}}}

\newcommand\al{\alpha}

\newcommand\ga{\gamma}
\newcommand\de{\delta}

\newcommand\la{\lambda}

\newcommand\si{\sigma}

\newcommand\om{\omega}

\newcommand\pa{\partial}

\renewcommand*{\Im}{\operatorname*{Im}}

\DeclareMathOperator{\supp}{supp}%
\newcommand\Ivar{\nu} 

\def\aside#1{}
\def\aside#1{\footnote{#1}}

\numberwithin{equation}{section}
\theoremstyle{plain}
\newtheorem{thm}{Theorem}[section]

\newtheoremstyle{remarks}{}{}{\rmfamily}{}{\bfseries}{}{8pt}
{\thmname{#1}\thmnumber{ #2}:\thmnote{ #3.}}%
\theoremstyle{remarks}

\newtheorem{rem*}{Remark}

\title{Multidimensional decay in van der Corput lemma}
\author{Michael Ruzhansky}
\address{ 
Department of Mathematics \endgraf
Imperial College London\endgraf
180 Queen's Gate \endgraf 
London SW7 2AZ \endgraf 
United Kingdom} 
\email{m.ruzhansky@imperial.ac.uk} 
\thanks{The author was supported by the Leverhulme
Research Fellowship and by EPSRC grant EP/E062873/01.
}
\date{\today}

\begin{document}
\maketitle
\begin{abstract}
In this paper we present a multidimensional version of
the van der Corput lemma where the decay of the oscillatory
integral is gained with respect to all space variables,
connecting the standard one-dimensional van der Corput lemma
with the stationary phase method.
\end{abstract}

\section{Introduction}

This paper is devoted to the estimates for the oscillatory
integrals of the type
\begin{equation*}\label{EQ:first}
I(\la)=\int_{\R^{N}}e^{i\la
\Phi(x)}a(x)\,dx\,,
\end{equation*}
where the support of $a\in C_0^\infty(\RN)$ is sufficiently
small. The estimate for $I(\la)$ as $\la\to\infty$ is well-known
in one dimension $N=1$. If $\Phi$ is real--valued and
$|\Phi^{(k)}(x)|\geq 1$ on the support of $a$, then the estimate
$|I(\lambda)|\leq c_k \lambda^{-1/k}$ holds when
$k\geq 2$, or when $k=1$ and $\Phi^\prime(x)$ is monotonic.
In this case the bound $c_k$ is also independent of
$\Phi$ and $\lambda$ (see e.g. Sogge \cite{Sog}
or Stein \cite{Ste}), and the decay rate is sharp.
This result plays a crucial role in various areas of analysis.
For example, it is closely related to sublevel set estimates
of the form
$$
\textrm{ meas } \{s\in \supp a: |\Phi(s)|\leq t\}\leq
 c_k t^{1/k},
$$
where $\Phi$ is a function as above, with numerous applications
in partial differential equations, 
microlocal analysis, harmonic analysis, etc.

A multidimensional version of these results would be of
great value, but presents many difficulties. 
It is known that for dimensions $N\geq 1$, if,
for example, $\partial^\alpha\Phi\geq 1$ on $\supp a$, then
$|I(\lambda)|\leq c_\alpha|\lambda|^{-1/|\alpha|}.$ The decay rate
here is sharp, but the constant $c_\alpha$ may depend 
on $\Phi$ and the estimate does not scale well. Again,
such estimate is closely related to the multilinear sublevel
set problem (see e.g. Phong, Stein and Sturm \cite{PSS}).
Parameter dependent sublevel set estimates were recently
used by Kamotski and Ruzhansky \cite{KR} 
in the analysis of elliptic and hyperbolic systems
with multiplicities, to yield Sobolev space estimates for
relevant classes of oscillatory integrals and for the
solutions.

Recently, Christ, Carbery and Wright \cite{Ccw} and 
Carbery and Wright \cite{Cw}, 
proposed versions of van der Corput lemma for functions
of several variables, in formulations where the constant 
in the estimate is independent of the phase function. This
aspect is of significant importance for applications allowing
to investigate various perturbation and other properties of
appearing integrals. However, the decay rate of the 
corresponding oscillatory integral there is 
again essentially one--dimensional
because the non-degeneracy of only one 
(higher-order) derivative is assumed.

At the same time, the decay rate exhibited in many problems of
interest is better than one-dimensional. If one compares this
with the case of
non-degenerate stationary points of $\Phi$, the stationary phase
method will readily yield the decay rate 
$|I(\lambda)|\leq C \lambda^{-N/2}$. However, if a 
stationary point degenerates, the situation becomes much more
delicate (see e.g. H\"ormander \cite[Chapter 7]{Horm}), 
and no good estimates are available in general.

The aim of the present paper is to bridge the gap between 
van der Corput lemma and estimates provided by
the stationary phase method. On one
hand, the standard van der Corput lemma works well for 
degeneracies of high orders but produces only one-dimensional
decay rate. On the other hand, the stationary phase method
produces the multidimensional decay rate, but does not work
well for degenerate stationary points. 

The result of this paper
will yield a multidimensional decay rate for degenerate 
stationary points. We will identify a class of 
functions, for which this can be achieved. These functions
will have certain convexity type properties.
It is clear that certain convexity conditions are necessary
to ensure the multidimensional decay rate. In fact, conditions
of the one-dimensional van der Corput lemma guarantee that
the function (or some derivative of the function) is convex
in one dimension. Thus, it is natural that
an analogue of convexity also 
appears in several dimensions to ensure that we gain 
one-dimensional decays in all directions. It is then a
question of putting all these rates together to yield the
full multidimensional decay, which will turn out to be $N$--times
better than the standard van der Corput estimate.

In what follows we will also allow phase function $\Phi$ to be
complex valued and to depend on an arbitrary set of parameters.
These two situations often happen in applications to
partial differential equations, in particular in the analysis
of solutions represented as oscillatory integrals,
leading to the dispersive and to the subsequent Strichartz
estimates. Thus, the complex phase corresponds to the fact that
characteristics of the analysed evolution equations 
may be complex
(see e.g. Tr\`eves \cite{T}).
At the same time, the dependence of the phase and of the amplitude
on parameters is also essential, and is related to
uniform estimates of \cite{Ccw}. Also, in applications
to the Strichartz estimates for hyperbolic equations of high
orders considered by Ruzhansky and Smith \cite{RS}, a parameter
is essential to encode the information on low order perturbations
of the equation, in order to establish the dispersive estimates for
solutions uniformly over such perturbations. At the same time,
in hyperbolic equations with time dependent coefficients
(e.g. considered by Matsuyama and Ruzhansky \cite{MR}),
the parameter encodes the information on the perturbations of
the limiting behaviour of coefficients, again allowing to obtain
dispersive estimates uniformly over such perturbations. We will
leave out these and other applications outside the scope of
this short paper.

We will use the standard multi-index
notation $\alpha=(\alpha_1,\ldots,\alpha_N)$, denote
its length by $|\alpha|=\alpha_1+\cdots+\alpha_N$ and partial
derivatives by $\partial^\alpha=\partial_{x_1}^{\alpha_1}\cdots
\partial_{x_N}^{\alpha_N}$.
We will also
use the standard convention to denote all constants by
letter $C$ although they may have different values on 
different occasions. 

\section{Multidimentional van der Corput lemma}

The following theorem is the main result that
establishes the multidimensional decay rate for a class of
oscillatory integrals. 

\begin{thm}\label{THM:oscintthm}
Consider the oscillatory integral
\begin{equation*}\label{EQ:genoscint}
I(\la,\Ivar)=\int_{\R^{N}}e^{i\la
\Phi(x,\Ivar)}a(x,\Ivar)\chi(x)\,dx\,,
\end{equation*}
where $N\geq 1$\textup{,} 
and $\Ivar$ is a parameter. 
Let $\ga\ge2$ be an integer.
Assume that
\begin{enumerate}[label=\textup{(A\arabic*)}]
\item\label{HYP:mainoscintgbdd} there exists a 
sufficiently small $\de>0$ such that 
$\chi\in C^\infty_0(B_{\delta/2}(0))$, where
$B_{\delta/2}(0)$ is the ball with radius 
${\delta/2}$ around $0$\textup{;}
\item\label{HYP:mainoscintImPhipos} $\Phi(x,\Ivar)$ is a
complex valued function such that $\Im \Phi(x,\Ivar)\ge0$ for all
$x\in \supp\chi$ and all parameters $\Ivar$\textup{;}
\item\label{HYP:mainoscintFconvexfn} for some fixed
$z\in\supp\chi$,  the function
\begin{equation*}
F(\rho,\om,\Ivar):=\Phi(z+\rho\om,\Ivar), \; |\om|=1,
\end{equation*}
satisfies the following conditions. Assume that
for each $\mu=(\omega,\nu)$, function $F(\cdot,\mu)$ is of class
$C^{\gamma+1}$ on $\supp\chi$, and
let us write its $\gamma^\text{th}$ order
Taylor expansion in $\rho$ at $0$ as
\begin{equation*}\label{EQ:Fformwithremainder}
F(\rho,\mu)=\sum_{j=0}^\gamma a_j(\mu)\rho^j + 
R_{\gamma+1}(\rho,\mu)\,,
\end{equation*}
where $R_{\gamma+1}$
is the remainder term. Assume that 
we have
\begin{enumerate}[leftmargin=*,label=\textup{(F\arabic*)}]
\item\label{ITEM:AssumpF1} $a_0(\mu)=a_1(\mu)=0$ for all
$\mu$\textup{;}
\item\label{ITEM:AssumpF2} there exists a constant $C>0$ such that
$\sum_{j=2}^\ga\abs{a_j(\mu)}\ge C$ for all
$\mu$\textup{;}
\item\label{ITEM:AssumpF3} for each
$\mu$\textup{,} $\abs{\pa_\rho F(\rho,\mu)}$ is
increasing in $\rho$ for $0<\rho<\de$\textup{;}
\item\label{ITEM:AssumpF4} for each $k\leq \gamma+1$,
$\pa_\rho^kF(\rho,\mu)$ is bounded uniformly in
$0<\rho<\de$ and $\mu$\textup{;}
\end{enumerate}

\item\label{HYP:mainoscintAderivsbdd} for each multi-index $\al$
of length $\abs{\al}\le \big[\frac{N}{\ga}\big]+1$\textup{,} there
exists a constant $C_\al>0$ such that 
$\abs{\pa_x^\al a(x,\Ivar)}\le
C_\al$ for all $x\in \supp\chi$ and all parameters $\Ivar$.
\end{enumerate}
Then there exists a constant $C=C_{N,\ga}>0$ such that
\begin{equation}\label{EQ:oscintbound}
\abs{I(\la,\Ivar)}\le C(1+\la)^{-\frac{N}{\ga}}\quad\text{for all
}\; \la\in[0,\infty)
\textrm{ and all parameters } \Ivar.
\end{equation}
\end{thm}

Theorem \ref{THM:oscintthm} obviously
includes the case where $a$ and $\Phi$ depend on
different sets of parameters. In this case we
may let $\nu$ run over the whole space of
parameters.

We also note that assumption (A3), 
or rather \ref{ITEM:AssumpF3},
can be view as an analogue
of a convexity assumption. Indeed,
if $F$ is real valued, then~\ref{ITEM:AssumpF3}
implies that the second
order derivative $\pa_\rho^2 F(\rho,\mu)$ 
does not change sign for $0<\rho<\de$,
because $\pa_\rho F(0,\mu)=0$ by \ref{ITEM:AssumpF1}. 
In turn, condition \ref{ITEM:AssumpF1} is not restrictive, 
since $a_0(\mu)$ can be
taken out of the integral, and non-zero $a_1(\mu)$ would
actually give a faster decay rate.

\begin{proof}
It is clear that~\eqref{EQ:oscintbound} holds for $0\le\la\leq 1$
since $\abs{I(\la,\Ivar)}$ is bounded for such $\la$,
in view of assumptions \ref{HYP:mainoscintgbdd},
\ref{HYP:mainoscintImPhipos} and \ref{HYP:mainoscintAderivsbdd}.
So, we may 
consider the case where $\la\geq 1$. 
Let $z\in\R^N$ be as in
\ref{HYP:mainoscintFconvexfn}, and set
$x=z+\rho\om$, where $\om\in \SNm$,
$\rho>0$. For $N=1$, we
use ${\mathbb S}^0=\set{-1,1}$.
Then we can write
\begin{equation*}
I(\la,\Ivar)= \int_{\SNm}\int_0^\infty e^{i\la
\Phi(z+\rho\om,\Ivar)}a(z+\rho\om,\Ivar)\chi(z+\rho\om)
\rho^{N-1}\,d\rho\, d\om\,.
\end{equation*}
It suffices to
prove~\eqref{EQ:oscintbound} for the inner integral.

Choose a function $\theta\in C_0^\infty([0,\infty))$, 
$0\le\theta(s)\le1$ for
all $s$, such that $\theta(s)$ 
is identically~$1$ for $0\le s\le\frac{1}{2}$ and is
identically zero for $s\ge1$. Then with our notation
$F(\rho,\om,\Ivar)=\Phi(z+\rho\om,\Ivar)$, we split the inner
integral into the sum of the two integrals
\begin{gather*}
I_1(\la,\Ivar,\om,z)=\int_0^\infty e^{i\la
F(\rho,\om,\Ivar)}a(z+\rho\om,\Ivar)\chi(z+\rho\om)
\theta(\la^{\frac{1}{\ga}}\rho)\rho^{N-1}\,d\rho\,, \\
I_2(\la,\Ivar,\om,z)=\int_0^\infty e^{i\la
F(\rho,\om,\Ivar)}a(z+\rho\om,\Ivar)\chi(z+\rho\om)
(1-\theta)(\la^{\frac{1}{\ga}}\rho)\rho^{N-1}\,d\rho\,.
\end{gather*}

Let us first estimate $I_1=I_1(\la,\Ivar,\om,z)$. Since
$\theta(\la^{\frac{1}{\ga}}\rho)=0$ for
$\la^{\frac{1}{\ga}}\rho\ge1$, changing variable
$\tau=\la^{\frac{1}{\ga}}\rho$, we have
\begin{equation*}
\abs{I_1} \le
C\int_0^\infty\theta(\la^{\frac{1}{\ga}}\rho)\rho^{N-1}\,d\rho
=C\int_0^\infty \tau^{N-1}\la^{-\frac{N-1}{\ga}}
\theta(\tau)\la^{-\frac{1}{\ga}}\,d\tau,
\end{equation*}
which yields the following estimate for $I_1$:
\begin{equation}\label{eq:estI1}
|I_1| \le C\la^{-\frac{N}{\ga}}\int_0^1 \tau^{N-1}
\,d\tau \leq C\la^{-\frac{N}{\ga}}.
\end{equation}

In order to estimate $I_2=I_2(\la,\Ivar,\om,z)$, let us first 
establish a useful estimate for functions $F$ satisfying
condition \ref{ITEM:AssumpF3}. We claim that
under condition \ref{HYP:mainoscintFconvexfn},
or rather under \ref{ITEM:AssumpF1}--\ref{ITEM:AssumpF4},
there exist constants $C,C_m>0$ such that we have estimates
\begin{gather}
\abs{\pa_\rho F(\rho,\mu)}\ge C\rho^{\ga-1}
\label{EQ:F'lowerbound}\\
\text{and }\abs{\pa_\rho^mF(\rho,\mu)}\le
C_m\rho^{1-m}\abs{\pa_\rho F(\rho,\mu)},\label{EQ:F^(m)upperbound}
\end{gather}
for all $0<\rho<\de$\textup{,} all parameters
$\mu$, and all $m\leq\gamma+1$.
First, note that for $0<\rho\le1$ and $m=\gamma+1$,
estimate \eqref{EQ:F^(m)upperbound} follows from
\eqref{EQ:F'lowerbound} and assumption 
\ref{ITEM:AssumpF4}. So we may only consider $m\leq \gamma$.

Now, assumption \ref{ITEM:AssumpF2} implies that
\begin{equation}\label{EQ:pidefnandbound}
\pi(\rho,\mu):=\sum_{j=2}^\ga j\abs{a_j(\mu)}\rho^{j-1}\ge
C\rho^{\ga-1}\,.
\end{equation}
Thus, in order to prove~\eqref{EQ:F'lowerbound}, it suffices to show
that
\begin{equation}\label{EQ:F'bddbelowbypi}
\abs{\pa_\rho F(\rho,\mu)}\ge C\pi(\rho,\mu)\quad\text{for all
}0<\rho<\de \textrm{ and all parameters } \mu.
\end{equation}
For $1\le m\le \ga$, we have, using~\ref{EQ:Fformwithremainder},
\begin{equation}\label{EQ:TaylorexpfordmF}
\pa_\rho^mF(\rho,\mu)=\sum_{k=0}^{\ga-m}
\frac{(k+m)!}{k!}a_{k+m}(\mu)\rho^{k}+ R_{m,\ga-m}(\rho,\mu)\,,
\end{equation}
where $R_{m,\ga-m}(\rho,\mu)= \int_0^\rho
\pa_s^{\ga+1}F(s,\mu)\frac{(\rho-s)^{\ga-m}}{(\ga-m)!}\,ds$ is
the remainder term of the $(\ga-m)^{\text{th}}$ Taylor expansion of
$\pa_\rho^mF(\rho,\mu)$. By~\ref{ITEM:AssumpF4} and
\eqref{EQ:pidefnandbound}, we get that
\begin{equation}\label{EQ:estremainderm<ga}
\abs{R_{m,\ga-m}(\rho,\mu)}\le C_{\ga,m}\rho^{\ga-m+1}\le
C_{\ga,m}\pi(\rho,\mu)\rho^{2-m}\quad\text{for }0<\rho<\de\,.
\end{equation}
Hence, for $0<\rho<\de$, we have
\begin{multline*}
\abs{\pa_\rho F(\rho,\mu)}=\absBig{\sum_{k=0}^{\ga-1}
(k+1)a_{k+1}(\mu)\rho^{k} +
R_{1,\ga-1}(\rho,\mu)}\\\ge\absBig{\sum_{j=2}^\ga
ja_j(\mu)\rho^{j-1}}-\absBig{ R_{1,\ga-1}(\rho,\mu)}
\ge\absBig{\sum_{j=2}^\ga ja_j(\mu)\rho^{j-1}}-
C_\ga\pi(\rho,\mu)\rho\,.
\end{multline*}
It follows now from assumptions \ref{ITEM:AssumpF1}
and \ref{ITEM:AssumpF3} that
\begin{align*}
\abs{\pa_\rho F(\rho,\mu)} &=\max_{0\le\si\le\rho}\abs{\pa_\rho
F(\si,\mu)}\\ \ge&  \max_{0\le\si\le\rho}\absBig{\sum_{j=2}^\ga
ja_j(\mu)\si^{j-1}}-
\max_{0\le\si\le\rho}C_\ga\pi(\si,\mu)\si\\
=&\max_{0\le\bar\si\le1}\absBig{\sum_{j=2}^\ga
ja_j(\mu)\rho^{j-1}\bar{\si}^{j-1}}- C_\ga\pi(\rho,\mu)\rho\,,
\end{align*}
since $\pi(\si,\mu)\si=\sum_{j=2}^\ga j\abs{a_j(\mu)}\si^{j}$
achieves its maximum on $0\le\si\le\rho$ at $\si=\rho$.
Noting that
\begin{equation*}
\max_{0\le\bar\si\le1}\absBig{\sum_{j=2}^\gamma z_j\bar{\si}^{j-1}}
\quad\text{and}\quad \sum_{j=2}^\gamma \abs{z_j}
\end{equation*}
are both norms on $\C^{\gamma-1}$ 
and, hence, are equivalent, we immediately get
\begin{align*}
\abs{\pa_\rho F(\rho,\mu)} \ge &C\sum_{j=2}^\ga
j\abs{a_j(\mu)}\rho^{j-1}-
C_\ga\pi(\rho,\mu)\rho \\
\ge& (C-C_\ga\de)\pi(\rho,\mu)\geq C \pi(\rho,\mu)\,,
\end{align*}
for some constants $C>0$,
if $\de$ is sufficiently small.
This completes the proof of~\eqref{EQ:F'bddbelowbypi}.

To prove~\eqref{EQ:F^(m)upperbound}, we will use
the representation \eqref{EQ:TaylorexpfordmF}.
Since $1\le m\le\ga$,  it follows from the definition of
$\pi(\rho,\mu)$ that
\begin{equation*}
\absBig{\sum_{k=0}^{\ga-m}\frac{(k+m)!}{k!}\,a_{k+m}
(\mu)\rho^{k}} \le C_m\pi(\rho,\mu)\rho^{1-m}\,,
\end{equation*}
which, together with~\eqref{EQ:estremainderm<ga} and
\eqref{EQ:F'bddbelowbypi}, yields
\begin{equation*}
\abs{\pa_\rho^mF(\rho,\mu)}\le C_{m,\de}\rho^{1-m}\abs{\pa_\rho
F(\rho,\mu)}\quad\text{for }0<\rho<\de.
\end{equation*}
This completes the proof of the claimed estimates
\eqref{EQ:F'lowerbound} and 
\eqref{EQ:F^(m)upperbound}.

Let us now come back to the estimate for $I_2$. Define
the operator 
$$L:=(i\la\pa_\rho F(\rho,\om,\Ivar))^{-1}\frac{\pa}{\pa
\rho}$$ which clearly satisfies the useful identity
$
L(e^{i\la F(\rho,\om,\Ivar)})=e^{i\la F(\rho,\om,\Ivar)}\,.
$
Denoting the adjoint of $L$ by $L^*$, we have, for each $l\in\NUO$,
\begin{equation*}
I_2=\int_0^\infty e^{i\la
F(\rho,\om,\Ivar)}(L^*)^l[a(z+\rho\om,\Ivar)\chi(z+\rho\om)
(1-\theta)(\la^{\frac{1}{\ga}}\rho)\rho^{N-1}]\,d\rho\,.
\end{equation*}
Now,
\begin{equation*}
(L^*)^l=\Big(\frac{i}{\la}\Big)^l\sum C_{s_1,\dots,s_p,p,r,l}
\frac{\pa_\rho^{s_1}F\cdots\pa_\rho^{s_p}F}{(\pa_\rho
F)^{l+p}}(\rho,\om,\Ivar)\frac{\pa^r}{\pa\rho ^r}\,,
\end{equation*}
where the sum is over all integers $s_1,\dots,s_p,p,r\ge0$ such that
$s_1+\dots+s_p+r-p=l$. From \eqref{EQ:F'lowerbound} and
\eqref{EQ:F^(m)upperbound} it follows that
\begin{equation*}
\absBig{\frac{\pa_\rho^{s_1}F\dots\pa_\rho^{s_p}F}{(\pa_\rho
F)^{l+p}}(\rho,\om,\Ivar)}\le
C\rho^{p-s_1-\dots-s_p-l\ga+l}=C\rho^{r-l\ga}\,.
\end{equation*}
Also, it is easy to see that for $r\le[\frac{N}{\ga}]+1$, we have
\begin{equation}\label{EQ:derivsofintegrandwrtrho}
\absBig{\frac{\pa^r}{\pa \rho^r}[a(z+\rho\om,\Ivar)\chi(z+\rho\om)
(1-\theta)(\la^{\frac{1}{\ga}}\rho)\rho^{N-1}]} \le C_{N}\rho^{N-1-r}
\widetilde{\chi}(\la,\rho)\,,
\end{equation}
where $\widetilde{\chi}(\la,\rho)$ is a smooth function in~$\rho$ which
is zero for $\la^{\frac{1}{\ga}}\rho<\frac{1}{2}$. 
Let us now take 
$l=[\frac{N}{\ga}]+1$, so that $N-l\ga<0$. Then we can
estimate
\begin{align*}
\abs{I_2}\le& C_{N}\la^{-l}\int_0^\infty\sum C_{s_1,\dots,s_p,p,r,l}
\;\rho^{r-l\ga}\, \rho^{N-1-r}\;\widetilde{\chi}(\la,\rho)\,d\rho\\
\le& C_{N}\la^{-l}\int_{\frac{1}{2}\la^{-\frac{1}{\ga}}}^\infty
\rho^{N-1-l\ga}\,d\rho = C_{N}\la^{-l}
\Big[\frac{\rho^{N-l\ga}}{N-l\ga}\Big]^{\infty}_{\frac{1}{2}
\la^{-\frac{1}{\ga}}} =C_{N,\ga}\la^{-\frac{N}{\ga}}.
\end{align*}
Combining this estimate with estimate \eqref{eq:estI1} for
$I_1$, we obtain the desired
estimate~\eqref{EQ:oscintbound}. 
This completes the proof of the theorem.
\end{proof}

We note that in the proof we showed that if function $F$
satisfies conditions 
\ref{ITEM:AssumpF1}--\ref{ITEM:AssumpF4}, it also satisfies
estimates \eqref{EQ:F'lowerbound} and \eqref{EQ:F^(m)upperbound}.
A version of this part of the argument was discussed by 
Sugimoto~\cite{Sug} for real valued
\emph{analytic} functions without dependence on~$\mu$, where
the analysis was based on the Cauchy's integral
formula for analytic functions (see also
Randol~\cite{Randol} and Beals~\cite{Beals}). 
The proof that we give for 
\eqref{EQ:F'lowerbound} and \eqref{EQ:F^(m)upperbound}
extends it to the generality required for Theorem
\ref{THM:oscintthm}. 

In fact, let us also briefly indicate a smooth version of
these estimates. Suppose that a function $F(\cdot,\mu)$
is smooth in the first variable, and that it satisfies 
conditions \ref{ITEM:AssumpF1}--\ref{ITEM:AssumpF3},
as well as condition \ref{ITEM:AssumpF4} for all $m\in\N$.
Then we claim that for sufficiently small $\de>0$, estimates 
\eqref{EQ:F'lowerbound} and \eqref{EQ:F^(m)upperbound}
are satisfied also for all $m\in\N$.

Indeed, we already proved estimate 
\eqref{EQ:F'lowerbound} and we also proved
\eqref{EQ:F^(m)upperbound} for $m\leq \gamma$. 
It remains to consider the case $m>\ga$.
Since $\ga+1-m\le0$, from~\ref{ITEM:AssumpF4} it 
trivially follows
that for $0<\rho<\de$ we have a stronger estimate
\begin{equation*}
\abs{\pa_\rho^mF(\rho,\mu)}\le C_m\le
C_{m,\de}\rho^{\ga+1-m}\le C_{m,\de}\rho^{2-m}\abs{\pa_\rho
F(\rho,\mu)},
\end{equation*}
where the last inequality is a consequence of
\eqref{EQ:F'lowerbound}.

\end{document}